\documentclass[12pt,thmsa]{article}
\usepackage{amsfonts}
\usepackage{amssymb}
\newtheorem{teo}{Theorem}[section]

\newtheorem{obs2}[teo]{Remark}

\newtheorem{tea}{Theorem}[subsection]

\newtheorem{no2}[teo]{Note}

\newtheorem{no3}[tea]{Note}

\newcommand{\Gal}{{\rm Gal}}
\newcommand{\Frob}{{\rm Frob }}
\newcommand{\trace}{{\rm trace}}
\newcommand{\mod}{{\rm mod}}

\newcommand{\Q}{\mathbb{Q}}

\newcommand{\F}{{\mathbb{F}}}

\newcommand{\cond}{{\rm cond}}

\newcommand{\rep} {\bar{\rho}_\ell}
\newcommand{\repp} {\bar{\rho}'_\ell}
\newcommand{\repe} {\bar{\rho}_{E, \ell}}
\begin{document}
\title{{\bf Existence of non-elliptic $\mod \; \ell$ Galois representations for every $\ell > 5$}}

\author{Luis Dieulefait
\\
Dept. d'{\`A}lgebra i Geometria, Universitat de Barcelona;\\
Gran Via de les Corts Catalanes 585;
08007 - Barcelona; Spain.\\
e-mail: ldieulefait@ub.edu\\
 }
\date{\empty}

\maketitle

\vskip -20mm

\begin{abstract}
For $\ell = 3$ and $5$ it is known that every odd, irreducible, $2$-dimensional representation of $\Gal(\bar{\Q}/\Q)$ with values in 
$\F_\ell$ and determinant equal to the cyclotomic character must ``come from" the $\ell$-torsion points of an elliptic curve defined over $\Q$. We prove, by giving concrete counter-examples, that this result is false for every prime $\ell >5$. 
\end{abstract}

\section{Examples for every $\ell>7$}
In [SBT] it is shown that for $\ell= 3$ and $5$ every odd, irreducible, $2$-dimensional Galois representation 
of $\Gal(\bar{\Q}/\Q)$ with values in 
$\F_\ell$ and determinant the cyclotomic character is ``elliptic", i.e., it agrees with the representation given by the action of      $\Gal(\bar{\Q}/\Q)$ on the $\ell$-torsion points of an elliptic curve defined over $\Q$. \\
In this note we will show that this is false for every prime $\ell >5$, i.e., that for every such prime there exists a Galois
 representation verifying the above properties but ``non-elliptic", i.e., not corresponding to the action of Galois on torsion points of any elliptic curve defined over $\Q$.  We will show this by giving concrete examples of non-elliptic representations.
 For any prime $\ell >7$, the example will be constructed starting from a weight $4$ classical modular form, corresponding to a rigid Calabi-Yau threefold. The case of $\ell = 7$ will be treated separately in the next section.\\

We consider the cuspidal modular form $f \in S_4(25)$ (i.e., of weight $4$, level $25$, and trivial nebentypus) which has all eigenvalues in $\mathbb{Z}$ and whose attached Galois representations $\rho_{f,\ell}$ agree (cf. [Y]) with the Galois representations on the third {\'e}tale cohomology groups of the Schoen rigid Calabi-Yau threefold. This threefold is obtained (after resolving the singularities) from:
$$ Y : X_0^5 + X_1^5 + X_2^5 + X_3^5 + X_4^5 - 5 X_0 X_1 X_2 X_3 X_4 = 0 \subseteq \mathbb{P}^4 $$ 
We list the first eigenvalues $a_p$ of $f$ (for $p \neq 5$):
$$a_2 = 1; a_3 = 7; a_7 =  6; a_{11} =  -43$$
From now on we will assume $\ell >5$.
For any prime $\ell$ let $\bar{\rho}_\ell := \bar{\rho}_{f,\ell}$ be the residual $\mod \; \ell$ representation corresponding to $\rho_{f,\ell}$: it is unramified outside $5 \ell$, its conductor or Serre's level 
(defined as the prime-to-$\ell$ part of its Artin conductor) divides $25$ and it has values in $\F_\ell$ and determinant 
$\chi^3$ ($\chi$ denotes the $\mod \; \ell$ cyclotomic character). For every prime
$p \nmid 5 \ell$ we have: $\trace (\bar{\rho}_\ell (\Frob \; p )) \equiv a_p \pmod{\ell}$. \\
Let us show that for any $\ell > 5$, $\rep$ is (absolutely) irreducible. As explained in [DM], since $\rho_\ell$ is attached to a rigid Calabi-Yau threefold, as long as $\ell > 4$ and $\ell$ is not $5$ (so that $\ell$ is a prime of good reduction), if 
$\rep$ is reducible it must hold: 
$$ \rep \cong \epsilon \oplus \epsilon^{-1} \chi^3 $$
where $\epsilon$ is a character unramified outside $5$ (the same description follows also from the fact that the representation is attached to a weight $4$ cuspform). Since
$$\cond(\epsilon) \cond(\epsilon^{-1}) = \cond(\epsilon)^2 = \cond(\rep) \mid  25$$
we have: $\cond(\epsilon) \mid 5$. In particular if $\ell \neq 11$ we have $\epsilon(11) = 1$, therefore:
$$ -43 = a_{11} \equiv \trace(\rep (\Frob \; 11)) \equiv 1 + 11^3 \pmod{\ell} $$
But no prime $\ell > 5, \ell \neq 11$ divides $11^3+1-43$, and this proves irreducibility of $\rep$ for every $\ell>5$ except $11$.\\
To show that $\bar{\rho}_{11}$ is also irreducible, observe that since it is an odd representation, irreducibility and absolute irreducibility are equivalent for it. Thus, it is enough to find a prime $p \nmid 55$ such that the reduction modulo $11$ of the characteristic polynomial $x^2- a_p x + p^3$ is irreducible. Equivalently, we need the discriminant $\Delta_p = a_p^2 - 4 p^3$
to be a non-square modulo $11$. For $p=2$ we have $\Delta_2 = -31 \equiv 2 \pmod{11}$, which is a non-square, and this gives the irreducibility of $\bar{\rho}_{11}$.\\
We define $\repp := \rep \otimes \chi^{(\ell-3)/2}$, for any $\ell >5$. It is also irreducible and odd, but the advantage is that its determinant is $\chi$.\\
We ask the following: Is there any elliptic curve $E$ defined over $\Q$ such that the Galois representation $\repe$ corresponding to  its 
$\ell$-torsion points gives $\repp$ for some $\ell$ ?\\
Let us show  that this can not happen for any $\ell > 7$. \\
Suppose the opposite. Then, since $\repp$ is unramified at $2$ and $2 \not\equiv 1 \pmod{\ell}$, if $\repp \cong \repe$ it is known
(cf. [C] and [R]) that $\rho_{E, \ell}$, the $\ell$-adic representation corresponding to the $\ell$-adic Tate module of $E$, must be unramified or semistable at $2$. If it is unramified at $2$, let us call
 $c_2$ the trace of $\rho_{E,\ell} (\Frob \; 2)$. Since $| c_2 | \leq 2 \sqrt{2}$, it should be $c_2 = 0 , \pm 1$ or $\pm 2$.\\
 Comparing the traces of $\repp$ and $\repe$ at $\Frob \; 2$ we obtain:
 $$ a_2 2^{(\ell - 3)/2} \equiv 0, \pm 1 , \pm 2 \pmod{\ell} \qquad (I)$$
 If $\rho_{E, \ell}$ is semistable at $2$, since $\repp$ is modular and $\rho_{E, \ell}$ is also modular (because all elliptic curves over $\Q$ are modular) then we obtain from $\repp \cong \repe$ by level raising (cf. [G]):
 $$\trace ( \repp (\Frob \; 2) ) \equiv \pm (2+1) \equiv \pm 3 \pmod{\ell} $$
 Thus: 
$$ a_2 2^{(\ell - 3)/2} \equiv \pm 3 \pmod{\ell} \qquad (II)$$
We conclude from (I) and (II) that if for some $\ell >5$, $\repp$ comes from an elliptic curve, it must hold (recall that $a_2=1$):
$$ 2^{(\ell-3)/2} \equiv 0, \pm 1, \pm 2, \pm 3  \pmod{\ell} $$
Thus: $2^{\ell-3} \equiv 1, 4, 9 \pmod{\ell}$.\\
Applying Fermat's little theorem this gives: 
$ 2^{-2} \equiv 1,4, 9 \pmod{\ell}$, and this is false for every prime $\ell > 7$. \\
Remark: It is natural that our result does not apply to $\ell = 7$ since independently of the  value of $a_2$, we would never get a contradiction for $\ell = 7$ because $0 ,\pm 1, \pm 2 , \pm 3$ cover all possible values modulo $7$.\\

We conclude that for any prime $\ell > 7$ the representation $\repp$ is non-elliptic. 

\section{The case $\ell = 7$}
We will consider the example of a $\mod \; 7$ representation attached to a weight $2$ cusp form $f$ such that the field $\Q_f$ generated by its eigenvalues is not $\Q$, $7$ is totally split in $\Q$, the representation is irreducible but it can not come from any elliptic curve for the following simple reason: the conductor of the representation is too large, compared with the universal bounds for conductors (cf. [Si]) of elliptic curves defined over $\Q$. Recall that the $p$-part of the conductor of any elliptic curve over $\Q$ must divide: $256$ if $p=2$, $243$ if $p=3$ and $p^2$ if $p >3$.\\

Concretely, we take the following example: Let $f \in S_2(512)$ be the cusp form with $\Q_f = \Q(\sqrt{2})$ and eigenvalues:
$$a_3 = \sqrt{2}, a_5 = -2 \sqrt{2}, a_7 = -4, a_{11} = \sqrt{2}, a_{13} = 2 \sqrt{2}, ...., a_{29} = 6 \sqrt{2} $$
(we obtain this values from the website [St]).\\
The corresponding $\mod \; 7$ representation $\bar{\rho}_7$ has values in $\F_7$ and it is irreducible because the discriminant $\Delta_{29}$
is a non-square modulo $7$.\\
The conductor of any of the representations $\rho_\lambda := \rho_{f,\lambda}$ in the family attached to $f$ ($\lambda \nmid 2$),
is equal to $512$, the level of $f$. Therefore (cf. [C], pag. 789) the conductor of $\bar{\rho}_7$ is also $512$.\\
Since the $2$-part of the conductor of any elliptic curve is at most $256$, this implies that $\bar{\rho}_7$ can not correspond to any elliptic curve. Thus $\bar{\rho}_7$, whose determinant is the cyclotomic character, is non-elliptic. \\

Remarks:\\
1) We have computed another example, using [St], with $f \in S_2(2560)$, with the same properties: $\bar{\rho}_7$ irreducible, valued in 
$\F_7$, but non-elliptic for the same reason. The field $\Q_f$ corresponds to a root of the polynomial $x^4- 316 x^2 + 8836$ (in [St] one can obtain a list of eigenvalues of $f$), it is a quadratic extension of $\Q(\sqrt{7})$ in which $\sqrt{7}$ splits.\\

2) Observe that from the ``bounds for conductors" in [Se], since $7 \not\equiv \pm 1 \pmod{9} $ and 
$7 \not\equiv \pm 1 \pmod{p} $ for any $p>3$, every odd, irreducible Galois representation valued in $\F_7$ must have the
 $p$-part of its conductor bounded with the same bound holding for elliptic curves, for any $p>2$. Thus, it is only by searching for representations with ``large $2$-part of the conductor" that one can obtain a representation valued in $\F_7$ not satisfying the universal bounds for conductors of elliptic curves.\\
 On the other hand,
 since $7 \equiv -1 \pmod{8}$, the bound for the $2$-part of conductors given in [Se] does not apply to the case of representations with values in $\F_7$. 

\section{Bibliography}

[C] Carayol, H.,  {\it Sur les repr{\'e}sentations galoisiennes modulo $\ell$
attach{\'e}es aux formes modulaires}, Duke Math. J. {\bf 59} (1989) 785-801
\\

[DM] Dieulefait, L.; Manoharmayum, J., {\it Modularity of rigid Calabi-Yau threefolds over $\Q$}, 
 in ``Calabi-Yau Varieties and Mirror Symmetry", Fields
Institute Communications, {\bf 38}, AMS (2003)
\\

[G] Ghate, E., {\it An introduction to congruences between modular forms}, in ``Current Trends in Number Theory", proceedings of the International Conference on Number Theory, Harish-Chandra Research Institute, Allahabad,November 2000; S. D. Adhikari, S. A. Katre, B. Ramakrishnan (Eds.), Hindustan Book Agency (2002)
\\

[R] Ribet, K., {\it Report on $\mod \; \ell$ representations of $\Gal(\bar{\Q}/ \Q)$}, in ``Motives" (Seattle, WA, 1991), Proc. Sympos. Pure Math., {\bf 55}, Part 2, AMS (1994)
\\

[SBT] Shepherd-Barron, N.; Taylor, R, {\it
Mod 2 and mod 5 icosahedral representations}, J.A.M.S. {\bf 10} (1997) 283-298\\ 

[Se] Serre, J-P., {\it Sur les repr{\'e}sentations modulaires de degr{\'e}
$2$ de $\Gal(\bar{\mathbb{Q}} / \mathbb{Q})$}, Duke Math. J. {\bf 54}
(1987) 179-230\\

[Si] Silverman, J., {\it Advanced topics in the arithmetic of elliptic curves},
Springer-Verlag (1994) \\

[St] Stein, W., {\it The Modular Forms Explorer}, available at:\\
http://modular.fas.harvard.edu/mfd/mfe/ \\

[Y] Yui, N., {\it Update on the modularity of Calabi-Yau varieties},  in
``Calabi-Yau Varieties and Mirror Symmetry", Fields
Institute Communications, {\bf 38}, AMS (2003)\\

\end{document}